\documentstyle[12pt]{article}
\newlength{\abstractwidth}
\renewcommand{\thefootnote}{\fnsymbol{footnote}}
\renewcommand{\thanks}[1]{\footnote{#1}} 
\newcommand{\starttext}{ \setcounter{footnote}{0}
\renewcommand{\thefootnote}{\arabic{footnote}}}

\newcommand{\be}{\begin{equation}}
\newcommand{\bea}{\begin{eqnarray}}
\newcommand{\eea}{\end{eqnarray}} \newcommand{\ee}{\end{equation}}
 
 \def\ba{\begin{eqnarray}}
\def\ea{\end{eqnarray}}

\def\K{{\cal K}}

\def\r{\rho}

\def\ra{\rightarrow}

\def\o{\omega}

\def\log{\,{\rm log}\,}
\def\exp{\,{\rm exp}\,}

\def\o{\omega}

\def\l{\lambda}

\def\o{\omega}

\def\r{\rho}

\def\na{\nabla}

\def\ve{\varepsilon}
\def\ge{\geq}
\def\le{\leq}

\def\ov{\overline}
\def\ti{\tilde}

\def\P{{\bf P}}

\def\i{\infty}
\def\I{\int}

\def\ra{\rightarrow}

\def\na{{\nabla}}

\def\K{{K\"ahler}}

 \def\v{\vskip .1in}

\def\[{{\bf [}}
\def\]{{\bf ]}}

\def\pl{\partial}

\begin{document}
\starttext \baselineskip=18pt \setcounter{footnote}{0}
\newtheorem{theorem}{Theorem}
\newtheorem{lemma}{Lemma}
\newtheorem{proposition}{Proposition}
\newtheorem{definition}{Definition}
\newtheorem{corollary}{Corollary}
\parindent=0in

\baselineskip=15pt \setcounter{equation}{0} \setcounter{footnote}{0}

\begin{center}
{\Large \bf THE K\"AHLER-RICCI FLOW WITH  }

\v

{\Large \bf POSITIVE BISECTIONAL CURVATURE\footnote{Research supported in part by
National Science Foundation grants  DMS-02-45371, DMS-06-04805, DMS-05-14003,
and DMS-05-04285.}}
\\
\bigskip
\bigskip

{\large D.H. Phong$^*$, Jian Song$^{**}$, Jacob Sturm$^\dagger$ and
Ben Weinkove$^\ddagger$} \\

\bigskip

\begin{abstract}
{\small }  We show that the K\"ahler-Ricci flow on a manifold with positive first Chern class converges to a K\"ahler-Einstein metric  assuming positive bisectional curvature and certain stability conditions. 
\end{abstract}

\end{center}

\section{Introduction}
\bigskip

Let $X$ be a compact K\"ahler manifold of complex dimension $n$ with $c_1(X)>0$.   The Frankel conjecture, proved by Mori \cite{Mr} and Siu-Yau \cite{SY}, states that if $X$ admits a K\"ahler metric of positive bisectional curvature then it is biholomorphic to ${\mathbf P}^n$.  There has been much interest in obtaining a proof of this using the K\"ahler-Ricci flow:
\bea \label{krf}
\frac{\partial}{\partial t} g_{\bar k j} = g_{\ov{k} j} - R_{\bar k j}.
\eea
 By a result of Goldberg-Kobayashi \cite{GK},
this amounts to solving the following well-known `folklore' problem:   without using
the existence of a K\"ahler-Einstein metric, show that if a
K\"ahler metric has positive bisectional curvature then the K\"ahler-Ricci flow deforms it to a
K\"ahler-Einstein metric.

\v

We mention now some work related to this problem.  The case $n=1$ was settled by Hamilton \cite{H1}, 
Chow \cite{Cho} (see also Chen-Lu-Tian \cite{CLT}).
Bando \cite{B} and Mok \cite{Mk} showed that, in every dimension, the positivity of the bisectional 
curvature is preserved along the K\"ahler-Ricci flow.  Chen-Tian \cite{CT} used the Moser-Trudinger inequalities \cite{T1, TZ1} (see also \cite{PSSW1}) to show that 
 if there exists a K\"ahler-Einstein
metric then, starting at a metric with positive bisectional curvature, the flow converges to it.  Perelman  later showed, without any  curvature conditions, that the  flow converges to a K\"ahler-Einstein metric when one exists, and this was extended to K\"ahler-Ricci solitons by Tian-Zhu \cite{P2, TZ2}.  Using an
injectivity radius estimate of Perelman \cite{P1}, Cao-Chen-Zhu \cite{CCZ} showed that if the bisectional curvature is nonnegative then the Riemann
curvature tensor is bounded along the flow. 
Chen \cite{C} showed, using the  Frankel conjecture together with the flow, that an irreducible K\"ahler
manifold with positive orthogonal bisectional curvature is biholomorphic to $\P^n$. 

\v

In \cite{PS3},
it was shown that the folklore problem can be reduced to establishing various stability
conditions.  In this paper we succeed in making further progress along these lines.  We consider the following three conditions:

\v

\quad (A) \ The Mabuchi K-energy is bounded below on $\pi c_1(X)$;

\v

\quad (A') The Futaki invariant of $X$ is zero;

\v

\quad (B) \ Let $J$ be the complex structure of $X$, viewed as a tensor.  
Then the $C^{\infty}$ closure of the orbit of $J$ under the diffeomorphism group of $X$ does not contain
any complex structure $J_{\infty}$ with the property that the space of holomorphic vector fields with
respect to $J_{\infty}$ has dimension strictly higher than the dimension of the space of holomorphic vector
fields with respect to $J$.

\v

Conditions (A) and (A') and their relations to stability have been studied intensely in the last two decades, and for the definitions we refer the reader to the literature (see \cite{PS1}, for example).  Condition (B) was introduced in \cite{PS3}.  
It was shown there that if the curvatures along the K\"ahler-Ricci flow are uniformly bounded, 
and if (A) and (B) hold then the K\"ahler-Ricci flow converges exponentially fast to a K\"ahler-Einstein
metric. Note that the Riemann curvature tensor is bounded along the flow if the bisectional
curvature is nonnegative or, in the case of two complex dimensions, if we have the weaker condition of nonnegative Ricci
curvature with traceless curvature operator 2-nonnegative \cite{PS2}.   

\v

Our first result is as follows:

\begin{theorem}  \label{theorem1} \  Suppose there exists a K\"ahler metric $g_0$ on $X$ with nonnegative bisectional curvature  which is positive at one point.   Assume condition (A) holds. Then the K\"ahler-Ricci flow starting at $g_0$ converges exponentially fast in $C^{\infty}$ to a K\"ahler-Einstein metric.
\end{theorem}

\v
Now, at least {\it a priori}, the algebraic condition (A') is much weaker than (A).   Here, we strengthen the result of \cite{PS3} by replacing  (A) by condition (A').

\begin{theorem} \label{theoremPS}  Suppose that the Riemann curvature tensor is uniformly bounded along the K\"ahler-Ricci flow and that conditions (A') and (B) hold.  Then the K\"ahler-Ricci flow converges exponentially fast in $C^{\infty}$ to a K\"ahler-Einstein metric.
\end{theorem}

If $n\le 2$ we have:

\begin{theorem} \label{theorem3}
Assume $X$ has complex dimension 1 or 2,  $g_0$ has nonnegative bisectional curvature and condition (A') holds.  Then the K\"ahler-Ricci flow starting at $g_0$ converges exponentially fast in $C^{\infty}$ to a K\"ahler-Einstein metric.
\end{theorem}

\v

This result for $n=1$ has already been established by different methods  as mentioned above.
Theorem \ref{theorem3}  now shows that the folklore problem in complex dimension 2 can be reduced to a condition on the finite dimensional space of holomorphic vector fields.

\v

We remark that there are already proofs of Theorems \ref{theorem1} and \ref{theorem3} which first show the existence of a K\"ahler-Einstein metric and then apply the results of  \cite{CT}, \cite{P2}.  Indeed, Chen \cite{C} proved Theorem \ref{theorem1} by showing that the bisectional curvature along the flow approaches that of the Fubini-Study metric, concluding  that the manifold is ${\mathbf P}^n$, and then applying  \cite{CT}.   
A proof of Theorem \ref{theorem3} can be obtained by combining \cite{P2}
 with the result  that, in complex dimension 2, the vanishing of the Futaki invariant implies the existence of a K\"ahler-Einstein metric \cite{T1}.
We note that our proofs use primarily flow methods and in particular avoid showing first the existence of a K\"ahler-Einstein metric.

\v

A key step in the proofs of Theorems 1, 2 and 3 is to obtain a uniform lower bound for the first positive eigenvalue $\lambda$ of the $\bar \partial^{\dagger} \bar \partial$ operator on $T^{1,0}$ vector fields.  The idea of considering this eigenvalue along the K\"ahler-Ricci flow was introduced in \cite{PS3} and examined further in \cite{PSSW2}.   In Section 2 we show that certain curvature conditions imply the desired bound for $\lambda$.  In Sections 3, 4 and 5, we give the proofs of Theorems \ref{theorem1}, \ref{theoremPS} and \ref{theorem3} respectively.  Finally, in Section 6 we describe how the Deligne pairing can be used to show that the Futaki invariant vanishes in the case $n=1$.

\setcounter{equation}{0}

\section{Lower bounds for the $\bar \partial$ operator}

For a solution $g(t)$ of the K\"ahler-Ricci flow (\ref{krf}), we define the Ricci potential $u$ by 
${d\over dt}g_{\bar k j}=g_{\bar k j}-R_{\bar k j}=\pl_j\pl_{\bar k}u,$ where we normalize $u$ by
imposing the condition
$\I_Xe^{-u}\o^n=\I_X\o^n$.   Here, $\omega = \frac{\sqrt{-1}}{2} g_{\ov{k} j} dz^j \wedge d\ov{z}^k \in \pi c_1(X)$ is the K\"ahler form of $g(t)$. 

\v

In the following, we will make use of the estimates of Perelman \cite{P2} (see \cite{ST}):
\begin{enumerate}
\item[(i)] For a uniform $C>0$, we have $\displaystyle{ \| u \|_{C^0} + \| \nabla u \|_{C^0} + \| R \|_{C^0} \le C}$.
\item[(ii)] Let $\rho >0$ be given.  Then for all $x \in X$ and all $r$ with $0< r \le \rho$ we have
\be
\int_{B_r(x)} \omega^n > C' r^{2n},
\ee
for a uniform constant $C'>0$, where $B_r(x)$ is the geodesic ball  of radius $r$ centered at $x$ with respect to $g(t)$.
\item[(iii)] The diameter of $(X,g(t))$ is uniformly bounded.
\end{enumerate}

Define two time dependent inner products 
on $T^{1,0}$ by

\be \langle V,W\rangle_u \ = \ \I_Xg_{\bar k j}
V^j\overline{W^k}e^{-u}\o^n
\ \ {\rm and} \ \ 
\langle V,W\rangle_0\ = \ \I_Xg_{\bar k j}
V^j\overline{W^k}\o^n.
\ee

Since $u$ is uniformly bounded  the corresponding norms 
$\| \cdot \|_u$ and $\| \cdot \|_{0}$ are equivalent.  Let $\tilde{\lambda}= \tilde{\lambda}(t)$ and $\lambda= \lambda(t)$ respectively be
the smallest positive eigenvalues of the operators $\tilde{L} = - g^{i \bar j} \nabla_i \nabla_{\bar{j}} +
g^{i \ov{j}} \nabla_i u \nabla_{\ov{j}}$ and $L =  - g^{i \bar j} \nabla_i \nabla_{\bar{j}}$ acting on
$T^{1,0}$ vector fields.  Denote by $\eta$ the space of holomorphic vector fields on $X$.  Then
$\tilde{\lambda}$ is the largest number satisfying
\be\label{tilde}
 \I_X |\na_{\bar i}V^k|^2e^{-u}\o^n
\ \geq \ 
\ti\l\I_X |V^k|^2e^{-u}\o^n
\ee
for all $V$ with the property: $\langle V,\xi \rangle_u =0$ for all $\xi\in \eta$.  
Similarly, $\lambda$ is the largest number satisfying
\be\label{notilde}
 \I_X |\na_{\bar i}V^k|^2\o^n
\ \geq \ 
\l\I_X |V^k|^2\o^n
\ee
for all $V$ with the property: $\langle V,\xi \rangle_0=0$ for all $\xi\in \eta$.  
The following lemma shows that $\tilde{\lambda}$ and $\lambda$ are uniformly equivalent.

\begin{lemma} \label{lemmaequiv} There exist uniform positive constants $A_1$ and $A_2$ such that
\be
A_1 \tilde{\lambda} \le \lambda \le A_2 \tilde{\lambda}.
\ee
\end{lemma}

{\it Proof of Lemma \ref{lemmaequiv}:} 
Let $V\in T^{1,0}$ be a smooth vector field
such that 
$\langle V,\xi\rangle_0=0\ \ 
{\rm for\ all}\ \  \xi\in\eta$.  Write
\be
V=W+\xi_0 {\rm \ \ with\ \ } \xi_0\in\eta
{\rm \ \ and\ \ } \langle W,\xi\rangle_u =0 
\ \ {\rm for} \ {\rm all} \ \ 
\xi\in\eta.
\ee
Then
\be 0\ =\ \langle V,\xi_0\rangle_0\ = \ 
\langle W,\xi_0\rangle_0+\langle \xi_0,\xi_0\rangle_0,
\ee
and the Cauchy-Schwarz inequality implies

\be \langle\xi_0,\xi_0\rangle_0^2\ \leq \ 
\langle W,W\rangle_0\langle \xi_0,\xi_0\rangle_0.
\ee
Hence 
there exist $c_1,c_2>0$ such that
\be
c_1\langle \xi_0,\xi_0\rangle_u \ \leq \
\langle\xi_0,\xi_0\rangle_0\ \leq \ 
\langle W,W\rangle_0\leq c_2\langle W,W\rangle_u.
\ee

Thus
$$ \I_X |\bar\na V|^2\o^n\geq c_3
\I_X |\bar\na V|^2e^{-u}\o^n
=c_3\I |\bar\na W|^2e^{-u}\o^n
\geq
c_3\tilde\l\I_X |W|^2e^{-u}\o^n \qquad \qquad
$$
\be 
= {c_3\tilde\l \over 2}\langle W,W\rangle_u
+{c_3\tilde\l\over 2}\langle W,W\rangle_u
\geq {c_3\tilde\l\over 2}\langle W,W\rangle_u +{c_3c_1\tilde\l\over 2c_2}\langle \xi_0,\xi_0\rangle_u \ \geq \ 
c_4 \tilde{\lambda} \langle V,V\rangle_0,
\ee
and it follows that $\lambda \ge c_4 \tilde{\lambda}$, giving the first inequality.  The second inequality follows similarly.  Q.E.D.

\bigskip

We recall some notions of positivity.  A tensor $T_{\ov{j} i \ov{l} k}$ is \emph{Griffiths} nonnegative  if
\be
T_{\ov{j} i \ov{l} k} \ov{V^j} V^i \ov{W^l} W^k \ge 0
\ee
for all vectors $V, W \in T^{1,0}$.  For brevity we write $T_{\ov{j} i \ov{l} k} \ge_{Gr} 0$.   The condition of nonnegative bisectional curvature means $R_{\ov{j} i \ov{l} k} \ge_{Gr} 0$.  We say that a tensor $T_{\ov{j} i \ov{l} k}$ is \emph{Nakano} nonnegative  if 
\be
T_{\ov{j} i \ov{l} k}  \ov{\zeta^{jl}} \zeta^{ik} \ge 0,
\ee
for all tensors $\zeta \in T^{1,0} \otimes T^{1,0}$, and we write $T_{\ov{j} i \ov{l} k} \ge_{Na} 0$ for short.

\bigskip

Next, we show that under a positive curvature condition, the eigenvalue $\lambda$ can be bounded below away from zero.

\begin{lemma} \label{lemmalowerbound} Suppose that a K\"ahler metric $g$ satisfies
\begin{equation} \label{eqnNa}
R_{\ov{j} i \ov{l} k} + R_{\ov{j} i} g_{\ov{l} k} - c \,  g_{\ov{j} i} g_{\ov{l} k} \ge_{Na} 0,
\end{equation}
for some constant $c>0$.  Then $\lambda \ge c$.
\end{lemma}
{\it Proof of Lemma \ref{lemmalowerbound}:} \ 
Recall the commutation formulae:
\begin{eqnarray}
(\nabla_i \nabla_{\ov{l}} - \nabla_{\ov{l}} \nabla_i) V^k & = & g^{k\ov{m}} R_{ \ov{l} i  \ov{m}p} V^p \\  
(\nabla_i \nabla_{\ov{l}} - \nabla_{\ov{l}} \nabla_i) a_{\ov{j}} & = & g^{m\ov{q}} R_{ \ov{l}i  \ov{j}m}  a_{\ov{q}},
\end{eqnarray}
for a $T^{1,0}$ vector field $V$ and a $(0,1)$ form $a$.  Let $V$ be an eigenvector of the operator $L$ with eigenvalue $\lambda$.  Then
\be
- g^{i \ov{j}} \nabla_i \nabla_{\ov{j}} V^k = \lambda V^k.
\ee
Apply $\nabla_{\ov{l}}$ to obtain
\be
- g^{i \ov{j}} \nabla_{\ov{l}} \nabla_i \nabla_{\ov{j}} V^k = \lambda \nabla_{\ov{l}} V^k.
\ee
Using the commutation formulae we have
\be
 - g^{i \ov{j}} \nabla_i \nabla_{\ov{l}} \nabla_{\ov{j}} V^k + g^{i \ov{j}} g^{k\ov{m}} R_{ \ov{l} i \ov{m}p} \nabla_{\ov{j}} V^p+ g^{i \ov{j}} g^{m\ov{q}} R_{\ov{l}i   \ov{j}m}  \nabla_{\ov{q}} V^k = \lambda \nabla_{\ov{l}} V^k.
 \ee
Multiply by  $g^{r\ov{l}} g_{\ov{t}k} \nabla_r  \ov{V^t}$ to obtain
\begin{eqnarray} \nonumber
&& - g^{r\ov{l}} g_{\ov{t}k} g^{i \ov{j}} \nabla_r  \ov{V^t} \nabla_i \nabla_{\ov{l}} \nabla_{\ov{j}} V^k +  g^{r\ov{l}}  g^{i \ov{j}}  R_{ \ov{l} i \ov{t}p}  \nabla_r  \ov{V^t} \nabla_{\ov{j}} V^p+    g^{r\ov{l}} g_{\ov{t}k}  g^{m\ov{q}}  R_{ \ov{l}m} \nabla_r  \ov{V^t}  \nabla_{\ov{q}} V^k \\
 && \mbox{} = \lambda g^{r\ov{l}} g_{\ov{t}k} \nabla_r  \ov{V^t} \nabla_{\ov{l}} V^k.
\end{eqnarray}
From (\ref{eqnNa}), after integrating by parts:
\be \lambda \int_X | \nabla_{\ov{i}} V^k |^2 \omega^n \ge c \int_X | \nabla_{\ov{i}} V^k |^2 \omega^n + \int_X | \nabla_{\ov{i}}  \nabla_{\ov{j}} V^k|^2 \omega^n,
\ee
and hence $\lambda \ge c$.  Q.E.D.

\bigskip

Next, we show, under a slightly different curvature assumption, that the eigenvalue $\tilde{\lambda}$ can be bounded below.

\begin{lemma} \label{lemmalambdatilde}Suppose that a K\"ahler metric $g$ satisfies 
\be
R_{ \ov{j} i  \ov{l} k} + (1-c)  g_{\ov{j} i} g_{\ov{l} k} \ge_{Na} 0,
\ee
for some constant $c>0$.  Then $\tilde{\lambda} \ge c$.
\end{lemma}
{\it Proof of Lemma \ref{lemmalambdatilde}:} \ Let $V$ be an eigenvector of $\tilde{L}$ with eigenvalue $\tilde{\lambda}$.  Then
\be- g^{i \ov{j}} \nabla_i \nabla_{\ov{j}} V^k + g^{i \ov{j}} \nabla_{\ov{j}} V^k \nabla_i u = \tilde{\lambda} V^k.
\ee
Applying $\nabla_{\ov{l}}$ as before, using the commutation formulae and the definition of $u$ we have
\begin{eqnarray} \nonumber
&& - g^{i \ov{j}} \nabla_i \nabla_{\ov{l}} \nabla_{\ov{j}} V^k + g^{i \ov{j}} g^{k\ov{m}} R_{ \ov{l} i  \ov{m}p} \nabla_{\ov{j}} V^p+ g^{i \ov{j}} g^{m\ov{q}} R_{ \ov{l} i  \ov{j}m}  \nabla_{\ov{q}} V^k \\
&& \mbox{} + g^{i \ov{j}} \nabla_{\ov{l}} \nabla_{\ov{j}} V^k \nabla_i u +  \nabla_{\ov{l}} V^k - g^{i \ov{j}} R_{ \ov{l}i} \nabla_{\ov{j}} V^k  = \tilde{\lambda} \nabla_{\ov{l}} V^k.
 \end{eqnarray}
Multiply by  $g^{r\ov{l}} g_{\ov{t}k} \nabla_r  \ov{V^t}$ to obtain
\begin{eqnarray} \nonumber
&& - g^{r\ov{l}} g_{\ov{t}k} g^{i \ov{j}} \nabla_r  \ov{V^t} \nabla_i \nabla_{\ov{l}} \nabla_{\ov{j}} V^k +  (R_{ \ov{j} i \ov{l}k} + g_{\ov{l}k} g_{ \ov{j}i}) \nabla^{\ov{j}}  \ov{V^l} \nabla^i V^k  \\
 && + g^{r\ov{l}} g_{\ov{t}k} g^{i \ov{j}} \nabla_r  \ov{V^t}  \nabla_{\ov{l}} \nabla_{\ov{j}} V^k \nabla_i u  \mbox{} = \tilde{\lambda} g^{r\ov{l}} g_{\ov{t}k} \nabla_r  \ov{V^t} \nabla_{\ov{l}} V^k.
\end{eqnarray}
Integrating against $e^{-u} \omega^n$ we obtain
\be \tilde{\lambda} \int_X | \nabla_{\ov{i}} V^k |^2 e^{-u} \omega^n \ge c \int_X | \nabla_{\ov{i}}  V^k |^2 e^{-u} \omega^n + \int_X | \nabla_{\ov{i}} \nabla_{\ov{j}} V^k|^2 e^{-u} \omega^n,\ee
and hence $\tilde{\lambda} \ge c$.  Q.E.D.

\setcounter{equation}{0}

\section{Proof of Theorem 1}

For the proof of Theorem \ref{theorem1},  we will need a number of lemmas.

\begin{lemma} \label{lemmaric}
Suppose the Mabuchi K-energy is bounded below on $\pi c_1(X)$ and the bisectional curvature of $g_0$ is nonnegative.  Then along the K\"ahler-Ricci flow
$$ \| R_{\ov{k} j} - g_{\ov{k} j} \|_{C^0} \rightarrow 0,$$
as $t \rightarrow \infty$.
\end{lemma}
{\it Proof of Lemma \ref{lemmaric}:} \ 
By the results of Bando \cite{B} and Mok \cite{Mk}, the nonnegativity of the bisectional curvature is preserved along the K\"ahler-Ricci flow.  It follows that the bisectional curvatures, and hence the full curvature tensor of $g=g(t)$ is uniformly bounded along the flow.  The covariant derivatives of the curvature are also uniformly bounded along the flow.  From \cite{PS3}, the lower boundedness of the Mabuchi K-energy implies
\be \label{eqnintR}
\int_X | R_{\ov{k} j} - g_{\ov{k} j} |^2 \omega^n = \int_X |R-n|^2 \omega^n  \rightarrow 0,
\ee
as $t \rightarrow \infty$.  Assume for a contradiction that there is a sequence of points $x_i$ and times $t_i \rightarrow \infty$ with $| R_{\ov{k} j} - g_{\ov{k} j} |(x_i, t_i) \ge \ve>0$.  Then by Perelman's non-collapsing result and the bound on the derivative of the Ricci curvature we obtain for  uniform constants $r>0$ and $c>0$,
\be
\int_{B_r(x_i)} | R_{\ov{k} j} - g_{\ov{k} j} |^2 \omega^n \ge c \, r^{2n},
\ee
at each time $t_i$.  This contradicts (\ref{eqnintR}).
 Q.E.D.

\bigskip

We will use the following result from \cite{C} (Theorem 1.5), which is proved using the maximum principle.

\begin{lemma} \label{lemmachen} Suppose there exist constants $c_0>0$ and $\nu>1/2$ such that the following holds.
There is a K\"ahler metric $g_0$ satisfying
\be
R_{ \ov{j} i  \ov{l} k}(g_0) - c_0 ( (g_0)_{\ov{j} i} (g_0)_{\ov{l} k} + (g_0)_{\ov{j} k} (g_0)_{\ov{l} i})  \ge_{Gr} 0,
\ee
 and the solution of the K\"ahler-Ricci flow $g=g(t)$ starting at $g_0$ satisfies
\be
R_{\ov{j} i} \ge \nu g_{\ov{j} i},
\ee
at all times.  Then, along the K\"ahler-Ricci flow, $g=g(t)$ satisfies
\be
R_{ \ov{j} i  \ov{l} k} - c_t (  g_{\ov{j} i} g_{\ov{l} k} + g_{\ov{j} k} g_{\ov{l} i}) \ge_{Gr} 0,
\ee
for $c_t>0$ with $\lim_{t \rightarrow \infty} c_t = (2\nu-1)/(n+1)>0$.
\end{lemma}

We will also need the following lemma:

\begin{lemma} \label{lemmademailly}
Suppose that the curvature of a K\"ahler metric $g$ satisfies
\be
R_{\ov{j} i \ov{l} k} - c  g_{\ov{j} i} g_{\ov{l} k}  \ge_{Gr} 0,
\ee
for some constant $c>0$.  Then
\begin{equation} \label{RR}
R_{\ov{j} i \ov{l} k} + R_{\ov{j} i} g_{\ov{l} k} - n c g_{\ov{j} i} g_{\ov{l} k}  \ge_{Na} 0.
\end{equation}
\end{lemma}
{\it Proof of Lemma \ref{lemmademailly}:}  This result is an application of the argument of \cite{D}, Proposition 10.14.  
It requires Lemma 10.15 from \cite{D}:

\begin{lemma} \label{lemmademailly2}
Let $q\ge 3$ be an integer and let $x^{\lambda}, y^{\lambda}$ for $1\le \lambda \le n$ be 
complex numbers.  Let $U^n_q$ be the set of $n$-tuples of $q$th roots of unity and define complex numbers
$$x'_{(\sigma)} = \sum_{\lambda=1}^n x^{\lambda} \ov{\sigma_{\lambda}}, \quad y'_{(\sigma)} = 
\sum_{\lambda=1}^n y^{\lambda} \ov{\sigma_{\lambda}}, \quad \textrm{for each } \sigma = (\sigma_1, \ldots,
\sigma_n) \in U^n_q.$$ Then for every pair $(\alpha, \beta)$ with $1 \le \alpha, \beta \le n$, the
following holds:
\begin{equation} \label{eqndemailly}
q^{-n} \sum_{\sigma \in U^n_q} x'_{(\sigma)} \ov{y'_{(\sigma)}} \sigma_{\alpha} \ov{\sigma_{\beta}} = 
\left\{ \begin{array}{ll} x^{\alpha} \ov{y^{\beta}}, \qquad & \textrm{if } \alpha \neq \beta \\
\sum_{\lambda =1}^n x^{\lambda} \ov{y^{\lambda}}, \quad & \textrm{if } \alpha=\beta.  \end{array} \right.
\end{equation}
\end{lemma}
{\it Proof of Lemma \ref{lemmademailly2}:} Although this lemma is already contained in \cite{D},
 we give the short proof here for the sake of completeness.  We only require the following elementary
claim: the coefficient of $x^{\lambda} \ov{y^{\mu}}$ in the left hand side of (\ref{eqndemailly}) is
$q^{-n} \sum_{\sigma \in U^n_q} \sigma_{\alpha} \ov{\sigma_{\beta}} \, \ov{\sigma_{\lambda}} \sigma_{\mu}$,
and this is equal to 1 if $\{ \alpha, \mu \} = \{ \beta, \lambda \}$ and  0 otherwise.  Indeed, for the
second alternative, assume without loss of generality that $\alpha \notin \{ \beta, \lambda \}$ and then
observe that
\begin{equation} \label{eqndemailly3}
\sum_{\sigma \in U^n_q} \sigma_{\alpha} \ov{\sigma_{\beta}} \, \ov{\sigma_{\lambda}} \sigma_{\mu} = \left\{ \begin{array}{ll} e^{2\pi i/q} \sum_{\sigma \in U^n_q} \sigma_{\alpha} \ov{\sigma_{\beta}} \, \ov{\sigma_{\lambda}} \sigma_{\mu}, \quad & \alpha \neq \mu \\
e^{4\pi i/q} \sum_{\sigma \in U^n_q} \sigma_{\alpha} \ov{\sigma_{\beta}} \, \ov{\sigma_{\lambda}} \sigma_{\mu}, \quad & \alpha = \mu. \end{array} \right.
\end{equation}
For (\ref{eqndemailly3}), replace $\sigma$ by the element of $U^n_q$ obtained by multiplying the $\alpha$ component of $\sigma$ by $e^{2\pi i/q}$.  Q.E.D.

\v

We may assume without loss of generality that we are calculating at a point where $g_{\ov{j} i} = \delta_{ji}$.  Fix $\zeta \in T^{1,0} \otimes T^{1,0}$.  We need to show 
\be
(R_{\ov{j} i \ov{l} k} + R_{\ov{j} i} g_{\ov{l} k} - nc g_{\ov{j} i} g_{\ov{l} k})  \ov{\zeta^{jl}} \zeta^{ik} \ge 0.
\ee
 Let $V_{(\sigma)}= V_{(\sigma)}^i \partial/\partial z^i$ be the vector with components $V_{(\sigma)}^i = \sum_{\lambda=1}^n \zeta^{i \lambda} \ov{\sigma_{\lambda}} \in \mathbf{C}$.  Let $W_{(\sigma)} = W_{(\sigma)}^k \partial/\partial z^k$ be the vector with components $W_{(\sigma)}^k = \sigma_k \in \mathbf{C}$.  Then, by assumption,
\begin{eqnarray} \nonumber
0 & \le & \sum_{i, j,k,l}  (R_{\ov{j} i \ov{l} k} - c g_{\ov{j} i} g_{\ov{l} k} )  q^{-n} \sum_{\sigma \in U^n_q}  \ov{V^j_{(\sigma)}} V^i_{(\sigma)}   \ov{W^l_{(\sigma)}} W^k_{(\sigma)} \\ \nonumber
& =& \sum_{i, j} \sum_{k\neq l}  R_{\ov{j} i \ov{l} k}   q^{-n} \sum_{\sigma \in U^n_q}  \ov{V^j_{(\sigma)}}  V^i_{(\sigma)}  \ov{\sigma_l} \sigma_{k}  + \sum_{i, j} \sum_{k=l} (R_{\ov{j} i \ov{l} k} - c g_{\ov{j} i} g_{\ov{l} k} )  q^{-n} \sum_{\sigma \in U^n_q}  \ov{V^j_{(\sigma)}} V^i_{(\sigma)}    \ov{\sigma_l} \sigma_{k} \\
& =& \sum_{i, j} \sum_{k\neq l}  R_{\ov{j} i \ov{l} k}  \ov{\zeta^{jl}} \zeta^{ik} + \sum_{i, j,k} (R_{\ov{j} i}  - nc g_{\ov{j} i}  )    \ov{\zeta^{jk}} \zeta^{ik},
\end{eqnarray}
where we have made use of Lemma \ref{lemmademailly2}.  Hence
\begin{eqnarray} \nonumber
\lefteqn{(R_{\ov{j} i \ov{l} k} + R_{\ov{j} i} g_{\ov{l} k} - nc g_{\ov{j} i} g_{\ov{l} k})  \ov{\zeta^{jl}} \zeta^{ik}
 } \\ &= & \sum_{k} \sum_{i,j}  R_{\ov{j} i \ov{k} k}  \ov{\zeta^{jk}} \zeta^{ik}  + \sum_{i,j} \sum_{k\neq l} R_{\ov{j} i \ov{l}k}  \ov{\zeta^{jl}} \zeta^{ik}  + \sum_{i,j,k} (R_{\ov{j} i} - nc g_{\ov{j} i})  \ov{\zeta^{jk}}  \zeta^{ik} \ge 0,
\end{eqnarray}
since the first term is nonnegative by the assumption. Q.E.D.

\bigskip

We can now prove Theorem \ref{theorem1}.

\bigskip

\noindent
{\it Proof of Theorem \ref{theorem1}:} \ If the initial metric has nonnegative bisectional curvature which is positive at one point then the bisectional curvature along the flow immediately becomes positive everywhere \cite{B, Mk}.   From Lemma \ref{lemmaric} we see  that for some $T>0$ and $\nu>1/2$ we have  $R_{\ov{j} i} \ge \nu g_{\ov{j} i}$ when $t\ge T$.  Without loss of generality then, we may assume that for $t\ge 0$ the metric has positive bisectional curvature and $R_{\ov{j} i} \ge \nu g_{\ov{j} i}$.  From Lemmas  \ref{lemmachen}, \ref{lemmademailly} and \ref{lemmalowerbound} we see that the eigenvalue $\lambda$ is uniformly bounded away from zero. Since the Mabuchi K-energy is bounded below it follows from Theorem 2 of \cite{PSSW2} (or, since the curvature is bounded,  the results of \cite{PS3}) that the K\"ahler-Ricci flow converges exponentially fast to a K\"ahler-Einstein metric. Q.E.D.

\setcounter{equation}{0}

\section{Proof of Theorem 2}

Before we give the proof of Theorem \ref{theoremPS}, 
we recall the definition of a \K-Ricci soliton.  We say that a metric $g$ with K\"ahler form
$\o\in \pi c_1(X)$ is a K\"ahler-Ricci soliton if 
\be g_{\bar k j}-R_{\bar k j}\ = \ \pl_j\pl_{\bar k} u
\ee
for a smooth function $u$ with $\bar\na\bar\na u=0$, or
in other words if $\na^ju$ is a  holomorphic vector field.
If $g$ is a \K-Ricci soliton then $g(t)=\Psi(t)^*g$ is a 
solution to the K\"ahler-Ricci flow, where $\Psi(t)$ is the 1-parameter
subgroup of holomorphic automorphisms generated by
the vector field $\textrm{Re}(\na^j u)$. Sometimes, by abuse of notation,
we also call $g(t)$ a \K-Ricci soliton.

\v

Now we recall from \cite{PS3} that for a solution $g(t)$ of the K\"ahler-Ricci flow, the function $Y(t)=\int_X | \nabla u|^2 \omega^n$ satisfies
\be\label{dotY}\dot Y(t) \ \leq \ -2\lambda(t) Y(t)-2\lambda(t) \textrm{Fut}(\pi_t (\na^ju))\ - \ Z(t),
\ee
where
\be Z(t)=\I_X|\na u|^2(R-n)\ + \ \I_X\na^ju\na^{\bar k}u(R_{\bar k j}-g_{\bar k j})\o^n,
\ee
and $\textrm{Fut}(\pi_t (\nabla^j u))$ is the Futaki invariant of the orthogonal projection $\pi_t$ with respect to $\langle \ , \ \rangle_0$ of the vector field $\nabla^j u$ to the space $\eta$ of holomorphic vector fields. We have the following lemma.

\begin{lemma} \label{lemmasoliton}
If $g(t)$ is a \K-Ricci soliton then $\dot Y(t)=Z(t)=0$ for all $t\geq 0$.
\end{lemma}
{\it Proof of Lemma \ref{lemmasoliton}:}  Since $Y$ is unchanged by automorphisms it follows that $\dot Y (t)=0$. Compute
$$\I_X\na^ju\na^{\bar k}u(R_{\bar k j}-g_{\bar k j})\o^n = \ 
-\I_X\na^ju\na^{\bar k}u(\pl_j\pl_{\bar k}u)\o^n
$$
\be =  
\I_X (\na_j\na^ju)(\na^{\bar k} u\na_{\bar k}u)\o^n =
\I_X (n-R)|\na u|^2\o^n,
\ee
and hence $Z(t)=0$.  
  Q.E.D.

\bigskip

We will make use of the following result. 

\v

\begin{theorem} \label{theoremPSSW}
Suppose condition (A') holds  and that along the K\"ahler-Ricci flow we have
$Y(t) \rightarrow 0$
as $t \rightarrow \infty$ and $\lambda(t) \ge c$ for some uniform constant $c>0$.  Then the K\"ahler-Ricci flow converges exponentially fast in $C^{\infty}$ to a K\"ahler-Einstein metric.
\end{theorem}
{\it Proof of Theorem \ref{theoremPSSW}:} This follows from the arguments of Lemma 5 and Lemma 6 of \cite{PSSW2}. Indeed, one can easily check that the argument of Lemma 5 of \cite{PSSW2} shows that under the assumptions of Theorem \ref{theoremPSSW},  $\| R(t)-n \|_{C^0}$ converges exponentially fast to zero.  Now apply Lemma 6 of \cite{PSSW2} which states that if $\int_0^{\infty} \| R(t) - n \|_{C^0} dt < \infty$ then the K\"ahler-Ricci flow converges exponentially fast in $C^{\infty}$ to a K\"ahler-Einstein metric. Q.E.D.

\v

We can now give the proof of Theorem \ref{theoremPS}.

\v

{\it Proof of Theorem \ref{theoremPS}:}  It is shown in \cite{PS3} that if the Riemann curvature tensor is uniformly bounded along the flow and condition (B) holds then there is a uniform lower bound of $\lambda(t)$ away from zero.
 If $Y(t) \rightarrow 0$ as $t \rightarrow \infty$ then the required result will follow immediately from Theorem \ref{theoremPSSW}.   We assume for a contradiction that there is a constant $\varepsilon>0$ and a sequence of times $t_j \rightarrow \infty$ such that
$Y(t_j) \ge \varepsilon$ for all $j$.

\v

Since we have uniformly bounded curvature, diameter and injectivity radius  along the flow we can apply Hamilton's compactness theorem \cite{H2} to obtain (after passing to a subsequence)
diffeomorphisms $F_j:\tilde{X} \ra X$ such that $F_j^*g(t_j+t)$ converges to
a solution $\tilde g(t)$ of the K\"ahler-Ricci flow on $\tilde X$ which is the
same manifold as $X$, but with possibly a different complex structure $\tilde J$ (see \cite{PS3}).  The convergence of the metrics and their derivatives is uniform on compact subsets of $\tilde{X} \times [0, \infty)$.
  Moreover, 
$\ti g$ is a K\"ahler-Ricci soliton.  

\v

This last assertion follows from a theorem in \cite{ST}, but for our particular case, we can give  here a direct argument for the convenience of the reader.  
Given a solution $g(t)$ of the K\"ahler-Ricci flow one can make a change of
variable $t= -\log(1-2s)$ and define a new metric $h=h(s)$ by
$h(s) = (1-2s) g(t(s))$.  Then $h$ satisfies, in real coordinates, $\frac{\partial}{\partial s} h_{ij} = -2 R_{ij}$ for $s \in [0, 1/2)$.  Perelman \cite{P1} showed that the functional
$$\mu(h, \tau) = \inf \{ (2\tau)^{-n} \int_X ( 2\tau (R + | \nabla f|^2) + f - 2n) e^{-f} \omega^n \ | \  (2\tau)^{-n} \int_X e^{-f} \omega^n = \int_X \omega^n \},$$
where the metric quantities are those of $h$,
satisfies $\frac{d}{ds}\mu (h(s), 1/2 - s) \ge 0$.   By Perelman's estimates for the scalar curvature and Ricci potential, $\mu$ is uniformly bounded from above.  Since $\mu$ is invariant under diffeomorphisms, it follows that the solution of the Ricci flow $\tilde{h}(s)$ corresponding to the limit solution $\tilde{g}(t)$ has $\mu (\tilde{h}(s), 1/2 - s)$ constant in $s$.   Hence (see for example \cite{KL}, section 12) $\tilde{h}$ satisfies $\tilde{R}_{ij} + \tilde{\nabla}_i \tilde{\nabla}_j f - \frac{1}{1-2s}\tilde{h}_{ij}=0$ for some $f=f(s)$ and it follows that $\tilde{g}$ is a K\"ahler-Ricci soliton, as required.  

\v

Now from (\ref{dotY}),
\be \dot Y(t_j+t)\ \leq \ -2\l Y(t_j+t) - \ Z(t_j+t).
\ee

Since $\lim_{j\to\i}Y(t_j+t)\ra \ti Y(t),  \lim_{j\to\i}Z(t_j+t)
=\tilde Z(t)$ uniformly for $t$ in any compact interval, we have

\be \dot{\tilde Y}(t)\leq -2\l \ti Y(t)\ - \ \tilde Z(t).
\ee

But Lemma \ref{lemmasoliton} says that $\dot{\tilde Y}(t)=\ti Z(t)=0$ so we get
$\tilde Y(t)=0$.   This contradicts the assumption that $Y(t_j) \ge \varepsilon$ for all $j$ and completes the proof of Theorem \ref{theoremPS}.  Q.E.D.

\setcounter{equation}{0}

\section{Proof of Theorem 3}

We now consider the case when $n\le 2$ and the Futaki invariant of $X$ vanishes.  The key lemma that makes use of one or two complex dimensions is as follows:

\begin{lemma} \label{lemmacurvature} Suppose $X$ has complex dimension $n \le 2$.   Then 
$$  R_{\ov{j} i \ov{l} k} \ge_{Gr} 0 \quad \Longleftrightarrow \quad R_{\ov{j} i \ov{l} k} \ge_{Na} 0.
$$
\end{lemma}
{\it Proof of Lemma \ref{lemmacurvature}:} \  The case $n=1$ is trivial.  Assume $n=2$ and that $g$ has nonnegative bisectional curvature.  We require $R_{\ov{j} i \ov{l} k}  \ov{\zeta^{jl}} \zeta^{ik} \ge 0$ for all $\zeta$.  Note that we only need  the inequality  for symmetric $\zeta$ since if we set $\nu^{ik} = (\zeta^{ik} + \zeta^{ki})/2$ then by the symmetry of the curvature tensor,
\be
R_{ \ov{j}i  \ov{l}k}  \ov{\zeta^{jl}} \zeta^{ik} = R_{ \ov{j} i \ov{l} k}  \ov{\nu^{jl}} \nu^{ik}.
\ee

We assume then that $\zeta$ is symmetric and of rank 2 (if $\zeta$ has rank 1 the result follows easily).  Make a linear change of complex coordinates so that $\zeta$ is the identity.  Denote these new coordinates by $z^1, z^2$.  Then
\be
R_{\ov{j}i   \ov{l}k}  \ov{\zeta^{jl}} \zeta^{ik} = R_{ \ov{1} 1 \ov{1}1} + R_{ \ov{2} 1 \ov{2}1} +  R_{ \ov{1}2  \ov{1}2} +  R_{ \ov{2}2  \ov{2}2}.
\ee
We will show that the right hand side is nonnegative.  Write $X = \partial/\partial z^1$ and $Y=  \partial/\partial z^2$.  Calculate
\begin{eqnarray} \nonumber 
0 & \le & R( \ov{X} - i \ov{Y},X+i Y, - i \ov{X} + \ov{Y}, i X + Y) \\ \nonumber
& = & R( \ov{X},X,  \ov{X},X)  + R( \ov{Y},X, \ov{Y}, X) +R( \ov{X},Y,\ov{X}, Y) + R( \ov{Y},Y, \ov{Y}, Y) \\
& =&   R_{ \ov{1}1  \ov{1}1} + R_{ \ov{2} 1 \ov{2}1} +  R_{ \ov{1}2  \ov{1}2} +  R_{ \ov{2}2  \ov{2}2},
\end{eqnarray}
where to go from the first to the second line, we have cancelled some terms using the symmetry of the curvature tensor.   Q.E.D.

\bigskip

\noindent
{\bf Remark} \ 
Note that positive bisectional curvature in dimension 2 is not equivalent to positive curvature in the Nakano sense.  Indeed, a K\"ahler manifold with $n\ge 2$ can never have positive Nakano curvature  because $R_{ \ov{j}i  \ov{l}k}  \ov{\zeta^{jl}} \zeta^{ik} =0$ for every skew-symmetric $\zeta$.

\bigskip

{\it Proof of Theorem \ref{theorem3}:} 
From Lemma \ref{lemmacurvature} and Lemma \ref{lemmalambdatilde} we see 
that $\tilde{\lambda}(t) \ge 1$ and so, by Lemma \ref{lemmaequiv}, $\lambda(t)$ is uniformly bounded below
away from zero along the flow.  We can now argue in the same way as in the proof of Theorem
\ref{theoremPS}. Q.E.D.

\setcounter{equation}{0}

\section{The Futaki invariant}

Suppose $X$ admits a  K\"ahler
metric of positive bisectional curvature. Then, by the 
Frankel conjecture, $X$ has a K\"ahler-Einstein metric
and hence the Futaki invariant of $X$ vanishes. 
In order to solve the folklore problem, one would
like to prove the vanishing of the Futaki invariant
without using the existence of a K\"ahler-Einstein metric.
In this section, we indicate how this can be done 
in the case $n=1$, using the Deligne pairing.

\v

\begin{proposition} \label{prop} If $n=1$, the Futaki invariant of $X$ vanishes.
\end{proposition}


{\it Proof of Proposition \ref{prop}.}  Let $K$ be the
canonical bundle of the Fano manifold $X$. We claim that the natural homomorphism
\be \textrm{Aut}^0(X)\ra \textrm{Aut}(\langle K,K\rangle)
\ee
is trivial, where $\langle \, , \, \rangle$ denotes the Deligne pairing.  Given this, the Futaki invariant vanishes by Theorem 1 of \cite{PS1}.

\v

We now prove the claim.  Let $V{\pl\over \pl z}$ be a holomorphic vector
field on $X$. The Poincar\'e-Hopf Theorem \cite{Mi} implies that $V$ has
two zeros (here we are using the topological classification of surfaces). Denote these zeros by $p$ and $q$, and assume for the moment
that $p\not=q$. Let $\Omega={1\over V}dz$. Then $\Omega$ is a meromorphic 1-form
on $X$ with simple poles at $p$ and $q$. After multiplying $V$ by a non-zero
scalar, we may assume that the residue at $p$ is $1$, and the residue at
$q$ is $-1$. Fix $z_0\in X$ with $z_0\not=p,q$ and let
\be f(z)\ = \exp{\left( \I_{z_0}^z \Omega \right)}.
\ee
Then $f$ is meromorphic, $f(p)=0$ and $f(q)=\i$. Moreover, if $\r_t$ is the
1-parameter family of biholomorphic maps generated by $\textrm{Im}(V)$, then
$f\circ \r_t=e^{-it/2}f$.

\v

Now let $\Omega_0= {1\over f}\Omega$ and $\Omega_1=f\Omega$. Then the divisor of $\Omega_0$ is
$-2p$ and the  divisor of $\Omega_1$ is $-2q$. In particular, the divisors
are disjoint so the Deligne pairing $\langle \Omega_0,\Omega_1\rangle$ is well-defined.
\v
Recall that if $f$ is a meromorphic function, and if $\Omega_0,\Omega_1$ are meromorphic
differential forms such that the divisors of $\Omega_0,\Omega_1$ and $f\Omega_1$ are pairwise
disjoint, then
\be \langle \Omega_0,f\Omega_1\rangle\ = \ f({\rm div}(\Omega_0))\langle \Omega_0,\Omega_1\rangle.
\ee
Thus we obtain

\be \langle \r_t^*\Omega_0,\r_t^*\Omega_1\rangle\ = \ 
\langle e^{it/2}\Omega_0,e^{-it/2}\Omega_1\rangle\ = \ 
\langle \Omega_0,\Omega_1\rangle.
\ee
\v

If $V$ has a double zero at the point $q$, then we let $p\in X$
with $p\not=q$, and define $h(z)=\I_{p}^z \Omega_1$ where $\Omega_1={1\over V}dz$.
Then $h$ vanishes at $p$ and has a simple pole at $q$ while
$\Omega_1$ has a double pole at $q$. Let $\Omega_0={1\over h^2}\Omega_1$. Then
$\Omega_0$ has a double pole at $p$. Let $\r_t$ be the
1-parameter family of biholomorphic maps generated by $\textrm{Im}(V)$.
Then 
\be \langle \r_t^*\Omega_0,\r_t^*\Omega_1\rangle\ = \ 
\langle F_t\Omega_0,\Omega_1\rangle,
\ee
for some meromorphic function $F_t(z)$. A simple
calculation shows that $F_t(q)=1$ and thus 
$\langle F_t\Omega_0,\Omega_1\rangle=\langle \Omega_0,\Omega_1\rangle$. Q.E.D.


\end{document}